\input amstex
\documentstyle{amsppt}
\magnification 1200
\NoBlackBoxes
\nopagenumbers
\topmatter
\title Vanishing theorems 
and syzygies for
 K3 surfaces and Fano varieties \endtitle
\author F. J. Gallego  and  B. P. Purnaprajna
\endauthor
\address{F.J. Gallego: Dpto. de \'Algebra,
 Facultad de Matem\'aticas,
 Universidad Complutense de Madrid, 28040 Madrid,
Spain}\endaddress
\email{gallego\@sunal1.mat.ucm.es}\endemail
\address{ B.P. Purnaprajna:
Dept. of Mathematics,
 Brandeis University, 
Waltham MA 02254-9110, USA}\endaddress
\email{purna\@max.math.brandeis.edu}\endemail
\thanks We are very pleased to thank David Eisenbud for his constant
encouragement, patience and advice. We are glad to thank N. Mohan
Kumar, Robert Lazarsfeld, Dale Cutkosky and Pradeep Shukla for their %%%%%%%%
constant encouragement and advice. It was Dale Cutkosky who pointed out to 
us that
our results for K3 surfaces should also go through for Fano varieties. 
\endthanks
\date May 26, 1996 \enddate

\abstract
In this article we prove some strong vanishing theorems on K3 surfaces. As
an application of them, we obtain higher syzygy results for K3 surfaces and
Fano varieties.
\endabstract

\endtopmatter
\document

\vskip .3 cm
 
\heading  1. Introduction\endheading

In this article we prove some vanishing theorems on K3 surfaces. An
application of the vanishing theorems is a result on higher syzygies for K3
surfaces and Fano varieties.

One part of our results fits a meta-principle stating that if $L$ is a line
bundle that is a product of $(p+1)$ ample and base point free line bundles
satisfying certain conditions, then $L$ satisfies the condition $N_p$ ( a
condition on the free resolution of the homogeneous coordinate ring of $X$
embedded by $L$). Other illustrations of this meta-principle have been given
in [GP1], [GP2] and [GP3]. The condition $N_p$ may be interpreted, through
Koszul cohomology, as a vanishing condition on a certain vector bundle.

The other part of our results provides strong vanishing theorems that imply,
in particular, the vanishing needed for $N_p.$ We also prove stronger
variants of the principle stated above for K3 surfaces and Fano varieties.

Before stating our results in detail, we recall some key results in this
area, namely the {\it normal generation }and {\it normal presentation }%
on K3 surfaces due to Mayer and St.Donat. Mayer and St. Donat proved that if 
$L$ is a globally generated line bundle on a K3 surface $X$ such that the
general member in the linear system is a non hyperelliptic curve of genus 
$%
g\geq 3,$ then $L$ is normally generated (in other words, the homogeneous
coordinate ring of $X$ in projective space $\bold{P}(H^0(L))$ is
projectively normal). If one further assumes that the general member in the
linear system is a non-trigonal curve which is not a plane quintic, then $L$
is normally presented (in other words, the homogeneous ideal $I_X$ defining 
$%
X$ in $\bold{P}(H^0(L))$ is generated by quadrics). Regarding results on
higher syzygies, it follows from a more general result of Ein and 
Lazarsfeld
[EL] that, if $L$ is a very ample bundle on a K3 surface then $L^{\otimes
(p+2)}$ satisfies $N_{p.}$ There are results on normal generation for Fano
threefolds due to Iskovskih. For details we refer the reader to [I].

In order to state our main theorem on syzygies we require the following:

Let $X$ be an irreducible projective variety and $L$ a very ample line
bundle on $X,$ whose complete linear series defines the morphism

$$
\Phi _L:X\longrightarrow \bold{P}(H^0(L)).
$$

Let $S=\bigoplus\limits_{m=0}^\infty S^mH^0(X,L)$ and $R(L)=\bigoplus%
\limits_{m=0}^\infty H^0(X,L^{\otimes m}).$ Since $R(L)$ is a finitely
generated module over $S,$ it has a minimal graded free resolution. We say
that the line bundle $L$ satisfies $N_p$, if $I_X$ is generated by quadrics
and the matrices in the free resolution of $R$ over $S$ have linear entries
until the $(p-1)${\it th} stage. Using this notation one says that a
normally generated line bundle satisfies $N_0$ and that a normally 
presented
line bundle satisfies $N_1$.

We prove that if $L$ is a globally generated bundle on a K3 surface $X$ such
that the general member in the linear system $C\in $ $\mid L\mid $ is a
non-hyperelliptic curve of genus $g\geq 4,$ then $L^{\otimes (p+1)}$
satisfies $N_{p.}$ If we further assume that the general member in the
linear system is a non-trigonal curve of genus $g\geq 5$ which is not a
plane quintic, then $L^{\otimes p\text{ }}$satisfies $N_{p.}$ We also show
that on a Fano threefold of index $1,$ Picard number $1$and very ample
anticanonical bundle of sectional genus $g\geq 4,$ $L^{\otimes (p+1)}$
satisfies $N_p$ for any ample line bundle $L$. If we further assume that the
general member in $\mid L\otimes \Cal{O}_S\mid $ is a non-trigonal curve
of genus $g\geq 5$ which is not a plane quintic curve (where $S$ is a
general, hence smooth, member of $\mid L\otimes \Cal{O}_S\mid )$, then $%
L^{\otimes p\text{ }}$satisfies $N_{p.}$ We generalize these results to Fano
varieties of dimension $n$ with index $(n-2).$ Our results generalize the
results of St. Donat and Iskovskih and improve the bound given by Ein and
Lazarsfeld for K3 surfaces and Fano varieties. Our results for K3 surfaces
do not assume that $L$ is ample.

\headline={\ifodd\pageno\rightheadline \else\leftheadline\fi}
\def\rightheadline{\tenrm\hfil \eightpoint VANISHING THEOREMS ON K3 
SURFACES
 \hfil\folio}
\def\leftheadline{\tenrm\folio\hfil \eightpoint F.J. GALLEGO \& B.P. 
PURNAPRAJNA \hfil}
%\voffset=2\baselineskip 

We refer the reader to Sections 4 and 6 for the statements of our vanishing
theorems on K3 surfaces and Fano varieties.

\heading 2. Normal  generation and normal presentation on K3 surfaces
\endheading

In this section we prove some vanishing theorems on K3 surfaces. As a
consequence of them we recover well known results on normal generation and
normal presentation due to Mayer and St. Donat. On the other hand the proofs
of these vanishing theorems will serve as a warm-up for the sequel, letting
us introduce part of the machinary and ideas used in Sections 4 and 6 to
prove results regarding higher syzygies. Moreover, Proposition 2.2 and 2.4
will be the first steps (and indeed, the key steps) of the inductive process
leading towards our higher syzygy results.

First of all, we will introduce the setting in which we will work and some
elementary facts about line bundles on K3 surfaces. Throughout this article
we will work over the field of complex numbers.

In Sections 2, 4 and 5, $X$ will be a smooth K3 surface and $L=\Cal{O}%
_X(C)$ will be a globally generated line bundle on $X$ (for example, if $%
C^2>0,$ it follows from [St.D] that $L$ is globally generated). Also, by
taking $C$ general in $\mid L\mid $, we can assume by Bertini's theorem that 
$C$ is smooth (See [St.D]). For a globally generated line bundle $G$ one 
has
the associated vector bundle $M_G$ given by the short exact sequence:

$$
0\longrightarrow M_G\longrightarrow H^0(G)\otimes \Cal{O}%
_X\longrightarrow G\longrightarrow 0\text{. }(2.1.1)
$$

We will constantly abuse the notation and write $C$ in place of $L=\Cal{O%
}_X(C)$ ( for example, we will write $M_C$ instead of $M_L,$ $pC$ instead 
of 
$L^{\otimes p})$. \vskip .2 truecm (2.1.2) An elementary and useful fact is
that $H^1(\Cal{O}_X(C))=0$ for any irreducible curve on $X.$ If, in
addition, the genus $g$ of $C$ is bigger than or equal to $2,$ $H^1(\Cal{%
O}_X(rC))=0$ for all $r\geq 1.$

A typical situation we will often encounter is the following: We have a
vector bundle $E$ and a line bundle $L$ and we want the multiplication map

$$
\varphi \: H^0(E)\otimes H^0(L^{\otimes r}) \longrightarrow %
H^0(E\otimes L^{\otimes r})
$$
to surject. We make the following

\proclaim {Remark 2.1}
Let $E$ be a coherent sheaf and let $L$ be a line bundle on a variety $X.$
Consider the following commutative diagram

$$
\matrix
H^0(E)\otimes H^0(L)^{\otimes r} & = & H^0(E)\otimes H^0(L)^{\otimes r} \\ 
\downarrow ^{\alpha _r} &  & \downarrow  \\ 
H^0(E\otimes L)\otimes H^0(L)^{\otimes r-1} &  & \vdots  \\ 
\downarrow ^{\alpha _{r-1}} &  & \downarrow  \\ 
H^0(E\otimes L^{\otimes 2})\otimes H^0(L)^{\otimes r-2} &  &  \\ 
\downarrow ^{\alpha _{r-2}} &  &  \\ 
&  &  \\ 
\downarrow ^{\alpha _3} &  & \downarrow  \\ 
H^0(E\otimes L^{\otimes r-2})\otimes H^0(L)^{\otimes 2} &  & H^0(E)\otimes
H^0(L^{\otimes r}) \\ 
\downarrow ^{\alpha _2} &  & \downarrow ^\psi  \\ 
H^0(E\otimes L^{\otimes r-1})\otimes H^0(L) & @>{\alpha _1}>> & 
H^0(E\otimes
L^{\otimes r})
\endmatrix
$$

The multiplication map $\psi $ is surjective if the maps $\alpha _1,\alpha
_2,...,\alpha _r$ are surjective.\medskip\ 
\endproclaim

We will use the following lemma very often:

\proclaim{Lemma 2.1}
Let $X$ be a regular variety ( i.e, a variety such that $H^1(\Cal{O}%
_X)=0).$ Let $E$ be a vector bundle on $X$ and $L$ $=\Cal{O}_X\left(
C\right) $ a globally generated line bundle such that $H^1(E\otimes
L^{-1})=0.$ If the multiplication map 
$$
H^0(E\otimes \Cal{O}_C)\otimes H^0(L\otimes \Cal{O}_C)@>
\alpha>>H^0(E\otimes L\otimes \Cal{O}_C)
$$
surjects, then the map 
$$
H^0(E)\otimes H^0(L)@>{\beta }>>H^0(E\otimes L)
$$
also surjects.
\endproclaim

{\bf Proof: }We have the sequence 

$$
0\longrightarrow \Cal{O}_X\longrightarrow \Cal{O}_X(C)%
\longrightarrow \Cal{O}_C\left( C\right) \longrightarrow 0.
$$
Taking the global sections of the short exact sequence above and then
tensoring with $H^0(E)$ yields the following commutative diagram: 
$$
\matrix
H^0(E)\otimes H^0(\Cal{O}_X) & \hookrightarrow & H^0(E)\otimes H^0(L) & 
\longrightarrow & H^0(E)\otimes H^0(L\otimes \Cal{O}_C)\longrightarrow 0
\\ 
\downarrow &  & \downarrow &  & \downarrow \\ 
H^0(E) & \hookrightarrow & H^0(E\otimes L) & \longrightarrow & H^0(E\otimes
L\otimes \Cal{O}_C)\ .
\endmatrix
$$

The vertical left hand side arrow is surjective for trivial reasons. The
vertical right hand side arrow is surjective because of the surjectivity 
of $%
\alpha $ and the vanishing of $H^1(E\otimes L^{-1})=0.$ The exactness at 
the
right of the top horizontal sequence follows from the vanishing of $H^1(%
\Cal{O}_X).$ The surjectivity of $\beta ,$which is the vertical middle
arrow, is then obtained by chasing the diagram. \ $\square$ \medskip\ 

\proclaim{Proposition 2.2.}
Let $X$ be a K3 surface and $L=\Cal{O}_X(C)$ be a line bundle on it.
Assume further that the general member $C\in \mid L\mid $ is a smooth
non-hyperelliptic curve of genus $g\geq 3.$ Then 
$$
H^1(M_{pC}\otimes rC)=0\text{ for all }r\text{ and }p\geq 1
$$
\endproclaim

{\bf Proof:} We consider the sequence (2.1.1) with $G=pC,$we tensor it
with $rC$ and take global sections. Since $H^1(rC)=0$ by (2.1.2), we obtain

$$
H^0(pC)\otimes H^0(rC)@>{\psi }>>H^0((p+r)C)%
\longrightarrow H^1(M_{pC}\otimes rC)\longrightarrow 0 \ ,
$$
so the vanishing of $H^1(M_{pC}\otimes rC)$ is equivalent to the
surjectivity of $\psi .$ To show the surjectivity of $\psi $ we use Remark
2.1. According to the remark we need to check the surjectivity of several
maps. Here we will only show the surjectivity of

$$
H^0(pC)\otimes H^0(C)@>{\alpha }>> H^0((p+1)C);
$$
(note that this map corresponds to $\alpha _r$ in Remark 2.1. Similar
arguments work for the rest of the maps $\alpha _{r-1},...\alpha _1).$ By
(2.1.2) and Lemma 2.1, it is enough to check that 
$$
H^0(pC\otimes \Cal{O}_C)\otimes H^0(C\otimes \Cal{O}%
_C)\longrightarrow H^0((p+1)C\otimes \Cal{O}_C)
$$
surjects. Since $C\otimes \Cal{O}_C$ is the canonical divisor $K_C$ on $%
C $ the required surjection follows from Noether's theorem.\ $\square$ 
\medskip\ 

\proclaim{Corollary 2.3}
(St. Donat, A. Mayer): Let $X$ be a K3 surface and $L=\Cal{O}_X(C)$ be a
line bundle on it. Assume further that the general member $C\in \mid L\mid 
$
is a smooth non-hyperelliptic curve of genus $g\geq 3$. Then $L=\Cal{O}%
_X(C)$ is normally generated.
\endproclaim

{\bf Proof: } From (2.1.1) it suffices to prove that $%
H^1(M_L\otimes L^{\otimes r})=0$ for all $r\geq 1$. This follows from the
above lemma by taking $p=1.$\ $\square $\medskip\ 

If we impose extra conditions on $C$ we obtain a stronger vanishing result:

\proclaim {Proposition 2.4}
Let $X$ be a K3 surface and $L=\Cal{O}_X(C)$ be the line bundle
corresponding to the divisor $C,$ where $C$ is a non-trigonal, non-plane
quintic curve of genus $g\geq 4.$ Then, 
$$
H^1(M_{kC}^{\otimes 2}\otimes rC)=0\text{ for all }k\text{ and }r\geq 1
$$
\endproclaim
\medskip 
We recover from this proposition the following result of St. Donat:

\proclaim{Corollary 2.5}
(St. Donat): Let $X$ be a K3 surface and $L$ a line bundle as in the above
proposition. Then, $L$ is normally presented.
\endproclaim

{\bf {}Proof: }By Corollary 2.3, $L$ is normally generated. Then by [GL1],
it is enough to prove that $H^1(\bigwedge^2M_L\otimes L^{\otimes r})=0$ for
all $r\geq 1$. Since we are working over the field of complex numbers, it 
is
enough to prove that $H^1(M_L^{\otimes 2}\otimes L^{\otimes r})=0$. This
follows from the Propositions 2.2 and 2.4 by taking $k=1.$ \ $\square $ %
\medskip\ 

Proposition 2.4 gives us a stronger vanishing result than the one needed to
prove Corollary 2.5, which is $H^1(M_C^{\otimes 2}\otimes rC)=0.$ We prove 
Proposition 2.4 in several steps, starting by proving the vanishing just
mentioned.

\proclaim{Lemma 2.6}
Let $X$ be a K3 surface and $L=\Cal{O}_X(C)$ be the line bundle
corresponding to the divisor $C,$ where $C$ is a non-trigonal curve of 
genus 
$g\geq 4,$ which is not isomorphic to a smooth plane quintic$.$ Then, 
$$
H^1(M_C^{\otimes 2}\otimes rC)=0\text{ for all }r\geq 1
$$
\endproclaim 

{\bf {}\ Proof:} We prove the lemma when $r=1.$ The proof when $r\geq 2$ is
less complicated and follows from similar arguments. Tensoring (2.1.1) 
with $%
M_C\otimes C$ and taking global sections yields 
$$
H^0(M_C\otimes C)\otimes H^0(C)@>{\varphi }>>%
H^0(M_C\otimes 2C)\longrightarrow H^1(M_C^{\otimes 2}\otimes C)
$$
$$
\longrightarrow H^1(M_C\otimes C)\otimes H^0(C).
$$
The group $H^1(M_C\otimes C)$ vanishes by Proposition 2.2. Therefore the
vanishing of $H^1(M_C^{\otimes 2}\otimes C)$ is equivalent to the
surjectivity of $\varphi .$ To show the surjectivity of $\varphi $ we use
Lemma 2.1. We need to see that $H^1(M_C)=0$ and that 
$$
H^0(M_C\otimes C\otimes \Cal{O}_C)\otimes H^0(C\otimes \Cal{O}_C)%
@>{\alpha }>>H^0(M_C\otimes 2C\otimes \Cal{O}_C)
$$
surjects. The former follows from the vanishing of $H^1(\Cal{O}_X).$ To
prove the later we consider the following sequence (see [GP2] for details): 
$$
0\longrightarrow \Cal{O}_C\longrightarrow M_C\otimes \Cal{O}%
_C\longrightarrow M_{C\otimes \Cal{O}_C}\longrightarrow 0\text{ }(2.6.1).
$$
If we tensor (2.6.1) by $C$ and take global sections, we obtain 
$$\displaylines {
0 \longrightarrow H^0(K_C)\longrightarrow H^0(M_C\otimes C\otimes \Cal{%
O}_C)@>{\beta }>>H^0(M_{K_C}\otimes
K_C)\longrightarrow 0 \cr
\ \longrightarrow H^1(K_C)\longrightarrow H^1(M_C\otimes C\otimes \Cal{%
O}_C)\longrightarrow H^1(M_{K_C}\otimes K_C)\longrightarrow 0.
}$$
The map $\beta $ is surjective because $h^1(M_C\otimes C\otimes \Cal{O}%
_C)=g+1$ and $h^1(M_{K_C}\otimes K_C)=g$, where $K_C=C\otimes \Cal{O}_C$
is the canonical bundle on $C.$ We can therefore write the following
commutative diagram (we denote $E=M_C\otimes C$ and $F=M_{K_C}\otimes 
K_C)$: 
$$
\matrix
0\longrightarrow & H^0(K_C)\otimes H^0(K_C) & \longrightarrow & 
H^0(E\otimes 
\Cal{O}_C)\otimes H^0(K_C)&@>{\beta }>> & 
H^0(F)\otimes H^0(K_C) \\ 
& \downarrow &  & \downarrow ^\alpha & &\downarrow ^\nu \\ 
0\longrightarrow & H^0(K_C^{\otimes 2}) & \longrightarrow & H^0(E\otimes
K_C)&\longrightarrow & H^0(F\otimes K_C) \ .
\endmatrix
$$
The left hand side vertical arrow is surjective by Noether's theorem and 
the
right hand side vertical arrow is surjective by Petri's theorem. Thus $%
\alpha $ also surjects as we wished. \ $\square$ \medskip\ 

We used Petri's theorem in the course of the last proof. We would like now
to give an alternative proof in the case of non-bielliptic curves without
using Petri's theorem, which fits in a more general context. It suggests
that we can get ``more'' than just the surjectivity we need. In fact using
the technique in our proof presented below and building upon the work of
[GL2], the second author and G. Pareschi prove that the canonical ring of a
curve $C$ satisfying the hypothesis in Lemma 2.6 is Koszul. We would like 
to
sketch the proof of the surjectivity of the map under consideration 
(readers
who are familiar with Koszul conditions in terms of surjectivity of global
sections of bundles will recognize that the surjection in the statement of
Lemma 2.7 is the first step needed to show that the canonical ring of the
curve under consideration is Koszul).%\medskip\ 

The alternative proof therefore follows the same steps as the previous one,
except that we show the surjectivity of $\nu $ in the following way:

\proclaim{Lemma 2.7}
Let $C$ be a smooth curve of genus $g\geq 6$ which is neither trigonal nor
bielliptic$.$ Assume further that it is not isomorphic to a smooth plane
quintic$.$ Then the multiplication map : 
$$
H^0(M_{K_C}\otimes K_C)\otimes H^0(K_C)\longrightarrow H^0(M_{K_C}\otimes
K_C^{\otimes 2})
$$
is surjective.
\endproclaim

{\bf {}\ Proof:} It follows from the hypothesis on\ $C$ and Mumford-Martens
theorems (cf. [ACGH]; see also [GL2]) that there exists a divisor $%
D=x_1+x_2+\cdots +x_{g-1}$ with $h^0(D)=2$ such that $D$ and $K_C(-D)$ are
globally generated.

Also by [GL2] we have the following exact sequences:\medskip\ \ 
$$
0\longrightarrow M_{K(-D)}\longrightarrow M_K\longrightarrow \Sigma
_D\longrightarrow 0\text{ (1)}
$$

$$
0\longrightarrow \Cal{O}_C(-x_{g-2}-x_{g-1})\longrightarrow \Sigma
_D\longrightarrow \bigoplus\limits_{i=1}^{g-3}\Cal{O}_C(-x_i)%
\longrightarrow 0\text{ (2)\ .}
$$
If we tensor (1) and (2) by $K_C$ and take global sections, we obtain 
$$
0\longrightarrow H^0(M_{K_C(-D)}\otimes K_C)\longrightarrow
H^0(M_{K_C}\otimes K_C)\longrightarrow H^0(\Sigma _D\otimes
K_C)\longrightarrow 0\text{ (3)}
$$

$$
0\longrightarrow H^0(K_C(-x_{g-2}-x_{g-1}))\longrightarrow H^0(K_C\otimes
\Sigma _D)\longrightarrow
H^0(\bigoplus\limits_{i=1}^{g-3}K_C(-x_i))\longrightarrow 0\text{ (4).}
$$
The exactness of (3) and (4) at the right hand side is deduced as follows:
the key point is to see that $H^0(K_C\otimes \Sigma _D)=g-2.$ For this, one
compares the bounds of $h^1(K_C\otimes \Sigma _D)$ obtained from the long
exact sequences of cohomology following (3) and (4), keeping into account
that $h^0(M_{K_C}\otimes K_C)=g$ by Noether's theorem and that $%
h^1(M_{K_C(-D)}\otimes K_C)$ $=2$ by the base-point-free pencil trick (cf.
[ACGH], Section 3). The rest is just a matter of adding and subtracting
dimensions of vector spaces.

Now denote $K_C\otimes \Sigma _D$ by $\Gamma $ and $M_{K_C(-D)}\otimes K_C$
by $P.$ If we tensor (3) by $H^0(K_C)$ and consider the obvious
multiplication maps, we obtain

$$
\matrix
0\longrightarrow H^0(P)\otimes H^0(K_C) & \longrightarrow H^0(F)\otimes
H^0(K_C)\longrightarrow & H^0(\Gamma )\otimes H^0(K_C)\longrightarrow 0 \\ 
\downarrow ^\delta & \downarrow ^\nu & \downarrow ^\epsilon \\ 
0\longrightarrow H^0(P\otimes K_C) & \longrightarrow H^0(F\otimes
K_C)\longrightarrow & H^0(\Gamma \otimes K_C)\longrightarrow 0 \ .
\endmatrix
$$
Therefore, to obtain the surjectivity of $\nu ,$ it suffices to check that
both $\delta $ and $\epsilon $ surject. Since $M_{K_C(-D)}=K_C(-D)^{*}$, 
the
map $\delta $ is in fact the multiplication map 
$$
H^0(D)\otimes H^0(K_C)\longrightarrow H^0(K_C\otimes D)
$$
which surjects by base-point-free pencil trick and the fact that $h^0(D)=2.$

To see the surjectivity of $\epsilon$ we argue as follows: we tensor (4) 
by $%
H^0(K_C)$ and considering the corresponding multiplication maps we obtain
the commutative diagram

$$
\matrix
0\longrightarrow H^0(T)\otimes H^0(K_C) & \longrightarrow H^0(\Gamma
)\otimes H^0(K_C)\longrightarrow & \bigoplus\limits_{i=1}^{g-3}(H^0(\Lambda
)\otimes H^0(K_C))\longrightarrow 0 \\ 
\downarrow & \downarrow & \downarrow \\ 
0\longrightarrow H^0(T\otimes K_C) & \longrightarrow H^0(\Gamma \otimes
K_C)\longrightarrow & \bigoplus\limits_{i=1}^{g-3}H^0(\Lambda \otimes
K_C)\longrightarrow 0
\endmatrix
$$
where $T=K_C(-x_{g-2}-x_{g-1})$and $\Lambda =K_C(-x_i)$. We want the
vertical arrows at the sides to be surjective. Since $C$ is
non-hyperelliptic, non-trigonal, non-bielliptic curve which is not
isomorphic to a plane quintic, Mumford-Martens says that the complete 
linear
series associated to $K_C(-x_i)$ and $K_C(-x_{g-2}-x_{g-1})$ give a
birational map from $C$ onto its image. The surjectivity of the vertical
arrows on the sides follows from the following result of Castelnuovo (see
[ACGH], page 151), for $D^{\prime }=K_C(-x_i),$$K_C(-x_{g-2}-x_{g-1})$
and $%
l=1:$ Let $\mid D^{\prime }\mid $ be a complete base point free linear
series of dimension $r=r(D^{\prime })\geq 3$ on a smooth curve $C,$ and
assume that
\pagebreak 

the mapping

$$
\Phi _{D^{^{\prime }}}:C\longrightarrow P^r
$$
is birational onto its image. Then the natural map

$$
Sym^lH^0(C,\Cal{O}(D^{\prime }))\otimes H^0(C,K_C)\longrightarrow
H^0(C,K_C(lD^{\prime }))
$$
is surjective for $l\geq 0.$ $\square$ \medskip\ 

Now we combine all these elements to give the proof of Proposition 2.4:
\medskip
{\bf {}Proof }(of Proposition 2.4). We will prove the result when $r=1$
(the arguments for $r\geq 2$ are similar). Thus we want to show that $%
H^1(M_{kC}^{\otimes 2}\otimes $ $C)=0$ for all $k$ $\geq 1$ . Since $%
H^1(M_{kC}\otimes $ $C)=0$ by Proposition 2.2, $H^1(M_{kC}^{\otimes
2}\otimes $ $C)$appears as the cokernel of 
$$
H^0(M_{kC}\otimes C)\otimes H^0(kC)@>{\varphi }>>%
H^0(M_{kC}\otimes (k+1)C)\ .
$$
To show that $\varphi $ is surjective we again use Remark 2.1. We will only
show here the surjectivity of one of the maps (in fact the first step,
corresponding to $\alpha _r)$ appearing in the statement of the remark,
namely $H^0(M_{kC}\otimes C)\otimes H^0(C)@>{\lambda }>>
H^0(M_{kC}\otimes
2C)$.

From (2.1.1$)$ it follows that the vanishing of $H^1(M_{kC}\otimes
M_C\otimes C)$ would imply the surjectivity of $\lambda .$ On the other
hand, since $H^1(kC)=0$ by (2.1.2), 
$$
H^1(M_{kC}\otimes M_C\otimes C)
$$
is also the cokernel of 
$$
H^0(kC)\otimes H^0(M_C\otimes C)@>{\mu }>>%
H^0(M_C\otimes (k+1)C)\ .
$$
Using again Remark 2.1, one concludes from (2.1.1) and Lemma 2.6 the
surjectivity of $\mu .$\medskip\ 

\heading 3. General lemmas \endheading

In this section we develop some homological tools necessary to prove the
vanishing theorems in the next section. We recall two lemmas from [GP2]. 
The
first is connected to the following problem: Consider two globally generated
line bundles $L_1$ and $L_2$. We would like to relate the vanishing of the
cohomology of $H^1(M_{L_1}^{\otimes p+1}\otimes L_2)$ to the vanishing of
the cohomology of a similar bundle on a divisor $Y$of $X,$obtained by
restricting $L_1$ and $L_2$ to $Y.$ One sees in the next section that Lemma
3.1 plays a crucial role in proving the vanishing theorems for K3 surfaces.
When one is working over an algebraic surface, the lemma transfers the
problem of computing cohomology of an unstable bundle to that of a
semistable bundle, which in general is easier to compute.\medskip\

\proclaim{Lemma 3.1}
Let $X$ be a projective variety, let $q$ be a non negative integer and let 
$%
F_i$ be a globally generated line bundle on $X$ for all $1\leq i\leq q+1.$
Let $Q$ be an effective line bundle on $X$ and let $C$ be a reduced and
irreducible member of $\mid Q\mid .$ Let $R$ be a line bundle on $X$ such
that

$(3.1.1)H^1(F_i\otimes Q^{*})=0$

$(3.1.2)H^1(R\otimes \Cal{O}_C)=0$

$(3.1.3)H^1(M_{(F_{i_1}\otimes \Cal{O}_C)}\otimes \cdots \otimes
M_{(F_{i_{q^{\prime }+1}}\otimes \Cal{O}_C)}\otimes R)=0$ for all $0\leq
q^{\prime }\leq q.$

Then, for all $-1\leq q^{^{\prime \prime }}\leq q$ and any subset $%
\{j_k\}\subseteq \{i\}$ with $\#\{j_k\}=q^{^{\prime \prime }}+1$ and for 
all 
$0\leq k^{^{\prime }}\leq q^{^{\prime \prime }}+1$ 
$$
H^1(M_{F_{j_1}}\otimes M_{F_{j_2}}\otimes \cdots \otimes
M_{F_{j_{k^{^{\prime }}}}}\otimes M_{(F_{i_{k^{^{\prime }}+1}}\otimes 
\Cal{O}_C)}\otimes \cdots \otimes M_{(F_{i_{q^{^{\prime \prime
}}+1}}\otimes \Cal{O}_C)}\otimes R\otimes \Cal{O}_C)=0\ .
$$
\endproclaim

\medskip 

The next lemma deals roughly with the following situation: Consider in
addition to $L_1$ and $L_2,$ two more ``positive'' line bundles 
$L_1^{\prime
}$ and $L_2^{\prime }$ (in the sense that $L_i^{*}\otimes L_i^{\prime }$ is
an effective line bundle). We would like to relate the vanishing of the
cohomology of $M_{L_1^{\prime }}^{\otimes p+1}\otimes L_2 $ and $%
M_{L_1}^{\otimes p+1}\otimes L_2^{\prime }$ to the vanishing of the
cohomology of $M_{L_1}^{\otimes p+1}\otimes L_2$. The usefulness of these
constructions is quite clear. For example, they give us a way to prove that
if a line bundle $L$ satisfies the property $N_p$, then so does the tensor
product of $L$ with certain effective line bundle.\medskip\ 

\proclaim{Lemma 3.2}
Let $X$ be a projective variety, let $q$ be a nonnegative integer and let 
$%
F_i$ be a globally generated line bundle on $X$ for all $1\leq i\leq q+1.$
Let $Q$ be an effective line bundle on $X$ and let $C$ be a reduced and
irreducible member of $\mid Q\mid .$ Let $R$ be a line bundle on $X$ such
that

$(3.2.1)H^1(F_i\otimes Q^{*})=0$

$(3.2.2)H^1(R\otimes Q\otimes \Cal{O}_C)=0$

$(3.2.3)H^1(M_{(F_{i_1}\otimes \Cal{O}_C)}\otimes \cdots \otimes
M_{(F_{i_{q^{^{\prime }}+1}}\otimes \Cal{O}_C)}\otimes R\otimes Q)=0$
for all $0\leq q^{^{\prime }}\leq q$

If $H^1(M_{F_1}\otimes M_{F_2}\otimes \cdots \otimes M_{F_{q+1}}\otimes R)
=0,
$ then

$$
H^1(M_{F_1}\otimes M_{F_2}\otimes \cdots \otimes M_{F_{q+1}}\otimes 
R\otimes
Q)=0.
$$
\medskip\ 
\endproclaim

We now prove a general lemma that will be used to prove vanishing theorems
for K3 surfaces in the next section.

\proclaim{Lemma 3.3}
Let $X$ be a projective variety and $L$ a globally generated line bundle.
Assume that

$(3.3.1)$ $H^1(M_L^{\otimes j}\otimes M_{L^{\otimes r}}^{\otimes i}\otimes
L^{\otimes p})=0$ for all $i,j$ such that $i+j=p$ and for all $r\geq 1.$

$(3.3.2)$ $H^1(M_L^{\otimes (p+1)}\otimes L^{\otimes (p+k)})=0$ for all $%
k\geq q_0\geq 0.$

Then $H^1(M_{L^{\otimes r}}^{\otimes i}\otimes M_L^{\otimes j}\otimes
L^{\otimes (p+k)})=0$ for all $i,j$ such that $i+j=p+1$, $k\geq q_0\geq 0$
and $r\geq 1.$
\endproclaim

{\bf {}\ Proof: }By induction on $i:$ If $i=0,$ we have $H^1(M_L^{\otimes
(p+1)}\otimes L^{\otimes (p+k)})=0$ by hypothesis. Assume the lemma to be
true for $i-1$ and we will prove it for $i.$ So we want 
$$
H^1(M_L^{\otimes j}\otimes M_{L^{\otimes r}}^{\otimes i}\otimes L^{\otimes
(p+k)})=0\text{ where }i+j=p+1\text{, }k\geq q_0\geq 0\text{ and }r\geq 1.
$$
By induction hypothesis, one has the vanishing $H^1(M_L^{\otimes
(j+1)}\otimes M_{L^{\otimes r}}^{\otimes (i-1)}\otimes L^{\otimes (p+k)})
=0$
for all $
% i+j=p+1$, $k\geq q_0\geq 0$ and $r\geq 1.$ But $H^1(M_L^{\otimes j}
\otimes
M_{L^{\otimes r}}^{\otimes i}\otimes L^{\otimes (p+k)})$ sits in the
following exact sequence obtained from (2.1.1):

$$
\displaylines {
H^0(M_L^{\otimes j}\otimes M_{L^{\otimes r}}^{\otimes (i-1)}\otimes
L^{(p+k)})\otimes H^0(L^{\otimes r})\longrightarrow H^0(M_L^{\otimes
j}\otimes M_{L^{\otimes r}}^{(i-1)}\otimes L^{(p+k+r)})\  \cr 
\longrightarrow H^1(M_L^{\otimes j}\otimes M_{L^{\otimes r}}^{\otimes
i}\otimes L^{\otimes (p+k)})\longrightarrow H^1(M_L^{\otimes j}\otimes
M_{L^{\otimes r}}^{\otimes (i-1)}\otimes L^{(p+k)})\otimes H^0(L^{\otimes
r}) \ .
}
$$
  
The last term of the above sequence is zero by hypothesis, since $i+j-1=p.$
In view of Remark 2.1, it is enough to prove that

$$
H^0(M_L^{\otimes j}\otimes M_{L^{\otimes r}}^{\otimes (i-1)}\otimes
L^{(p+k)})\otimes H^0(L)^{\otimes r}\longrightarrow H^0(M_L^{\otimes
j}\otimes M_{L^{\otimes r}}^{(i-1)}\otimes L^{(p+k+r)})\longrightarrow 0.
$$
We will show the first step, namely that 
$$
H^0(M_L^{\otimes j}\otimes M_{L^{\otimes r}}^{\otimes (i-1)}\otimes
L^{(p+k)})\otimes H^0(L)\longrightarrow H^0(M_L^{\otimes j}\otimes
M_{L^{\otimes r}}^{(i-1)}\otimes L^{(p+1+k)})
$$
surjects; the others are similar. Observe that the cokernel of the above
multiplication map is 
$$
H^1(M_L^{\otimes j+1}\otimes M_{L^{\otimes r}}^{\otimes (i-1)}\otimes
L^{(p+k)}),
$$
which is zero by induction assumption. $\square$

\medskip\ 

\proclaim{Lemma 3.4}
Let $X$ be a variety such that $H^1(\Cal{O}_X)=0.$ Let $L=\Cal{O}%
_X(C)$ be a globally generated line bundle. If

(1) $H^1(M_L^{\otimes p}\otimes L^{\otimes (p+k)})=0$

(2) $H^1(M_L^{\otimes p}\otimes L^{\otimes (p-1+k)})=0$

(3) $H^1(M_L^{\otimes p}\otimes M_{L\otimes \Cal{O}_C}\otimes L^{\otimes
(p+k)}\otimes \Cal{O}_C)=0$

then $H^1(M_L^{\otimes (p+1)}\otimes L^{\otimes (p+k)})=0.$
\endproclaim

{\bf {}\ Proof: }Observe that $H^1(M_L^{\otimes (p+1)}\otimes L^{\otimes
(p+k)})$ sits in the following exact sequence:

$$
\displaylines{
H^0(M_L^{\otimes p}\otimes L^{\otimes (p+k)})\otimes H^0(L)@>{\varphi 
}>> H^0(M_L^{\otimes p}\otimes L^{\otimes
(p+k+1)})\longrightarrow H^1(M_L^{\otimes (p+1)}\otimes L^{\otimes (p+k)})
\cr
\longrightarrow H^1(M_L^{\otimes p}\otimes L^{\otimes p})\otimes H^0(L) \ .
}
$$
Since $H^1(M_L^{\otimes p}\otimes L^{\otimes p})=0$ by hypothesis, it is
enough to prove that the multiplication map $\varphi $ surjects. In view of
Lemma 2.1, it is enough to check the vanishings of $H^1(M_L^{\otimes
p}\otimes L^{\otimes (p+k-1)})$ and $H^1(M_L^{\otimes p}\otimes M_{L\otimes 
\Cal{O}_C}\otimes L^{\otimes (p+k)}\otimes \Cal{O}_C)$, which
follows from hypothesis.
\ $\square$
\ \medskip\ 

\heading 4. Vanishings theorems on K3 surfaces \endheading

In this section we prove some strong vanishing theorems, namely Theorems 4.1
and 4.6, regarding certain vector bundles associated to a globally 
generated
line bundle on a K3 surface. As an application of Theorem 4.1, we obtain
higher syzygy results for K3 surfaces.\medskip\ 

\proclaim{Theorem 4.1}
Let $X$ be a K3 surface and $L$ a line bundle on $X$. Assume that the
general member $C\in \mid L\mid $ is a smooth non-trigonal, non-plane
quintic curve of genus $g\geq 5.$ Then 
$$
H^1(M_C^{\otimes i}\otimes M_{rC}^{\otimes j}\otimes (p+k)C)=0\text{ 
for all 
}i+j=p+1,\text{ }k\geq 0\text{ and }r\geq 1.\text{ }
$$ 
\endproclaim
\medskip\
To prove this theorem we need a number of intermediate results.

\proclaim{Lemma 4.2}
Let $X$ be a K3 surface and $L=\Cal{O}_X(C)$ be a line bundle such that
the smooth general member $C\in \mid L\mid $ is a non-hyperelliptic curve of
genus $g>3.$ Then the cohomology group $H^1(M_C^{\otimes p}\otimes 
(p+k)C)=0$
for all $k\geq 0$ and $p\geq 1.$
\endproclaim

{\bf {}\ Proof: }The proof is by induction on $p.$ If $p=1$ the vanishing
holds by Proposition 2.2. Assume it is true for $p-1$; then we have $%
H^1(M_C^{\otimes p-1}\otimes (p-1+k)C)=0$ for all $k\geq 0.$ We want to show
that $H^1(M_C^{\otimes p}\otimes (p+k)C)=0.$ Tensoring (2.1.1) (substitute 
$%
G=C)$ with $M_C^{\otimes (p-1)}\otimes (p+k)C$ and taking global sections 
yields
the long exact sequence:

$$\displaylines{
H^0(M_C^{\otimes p-1}\otimes (p+k)C)\otimes H^0(C)@>{\mu }>> H^0
(M_C^{\otimes
p-1}\otimes (p+1+k)C)\longrightarrow
\cr
\longrightarrow H^1(M_C^{\otimes p}\otimes (p+k)C)\longrightarrow
H^1(M_C^{\otimes p-1}\otimes (p+k)C)\otimes H^0(C).}
$$
The last term is zero by induction assumption. So it is enough to prove that
the multiplication map $\mu $ surjects.

Since $H^1(M_C^{\otimes p-1}\otimes (p-1+k)C)=0$ by induction, we may use
Lemma 2.1 to reduce the problem of seeing the surjectivity of the map 
$\mu $%
, which is a multiplication map of global sections of vector bundles on a K3
surface, to the problem of checking the surjectivity of a multiplication map
of sections of bundles on a curve. According to Lemma 2.1, it suffices to
show that

$$
H^0(M_C^{\otimes p-1}\otimes (p+k)C\otimes \Cal{O}_C)\otimes
H^0(C\otimes \Cal{O}_C)\longrightarrow H^0(M_C^{\otimes p-1}\otimes
(p+1+k)C\otimes \Cal{O}_C)
$$
surjects. This will follow from the vanishing 
$$
H^1(C,M_C^{\otimes p-1}\otimes \Cal{O}_C\otimes M_{K_C}\otimes
(p+k)K_C)=0,
$$
where $C\otimes \Cal{O}_C=K_C$ is the canonical bundle on $C.$ To see
the vanishing, we use Lemma 3.1. The first condition required by Lemma 3.1
is the vanishing of $H^1($ $\Cal{O}_X).$ The second is the vanishing of $%
H^1((p+k)K_C),$ which occurs because $p\geq 2.$ Finally, to verify the third
condition one has to check the vanishing of 
$$
H^1(M_{K_C}^{\otimes p-1}\otimes M_{K_C}\otimes (p+k)K_C).
$$

Since the bundle $M_{K_C\text{ }}$is stable (see [PR]) and the tensor
product of semistable bundles is semistable (see [Mi]), $M_{C\otimes 
\Cal{O}_C}^{\otimes p}\otimes (p+k)C\otimes \Cal{O}_C$ is
semistable. Therefore, it is enough to check that the slope $\mu
(M_{C\otimes \Cal{O}_C}^{\otimes (p)}\otimes (p+k)C\otimes \Cal{O}%
_C)>2g-2,$ where $g$ is the genus of $C$ and $\mu (E)=\frac{\deg E}{\text{%
rank} E}$ for a vector bundle $E.$ We have 
$$
\mu (M_{K_C}^{\otimes p}\otimes (p+k)K_C)=p\mu (M_{K_C})+(p+k)\deg (K_C),
$$
so we need $p(-2)+(p-1+k)(2g-2)>0,$ which is true for $p\geq 2$ and $g >3$. 
\
$\square$%
\medskip\ 

\proclaim{Lemma 4.3}
Let $X$ be a K3 surface and let $L=\Cal{O}_X(C)$ be a globally generated
line bundle such that the smooth general member $C\in \mid L\mid $ is a
non-hyperelliptic curve of genus $g>3.$ Then the cohomology group $%
H^1(M_C^{\otimes i}\otimes M_{rC}^{\otimes j}\otimes (p+k)C)=0$ for all $%
i+j=p$ and for all $k\geq 0$ and $r\geq 1.$
\endproclaim

{\bf {}Proof: }The proof is by induction on $i+j.$ If $i+j=1$ we want $%
H^1(M_C\otimes (p+k)C)$ and $H^1(M_{rC}\otimes (p+k)C)$ to vanish$,$ which
is true by Proposition 2.2. Let us assume that the result is true for $p-1.$
So we have $H^1(M_C^{\otimes i}\otimes M_{rC}^{\otimes j}\otimes 
(p-1+k)C)=0$
for $i+j= p-1.$ We want to prove the result for $p.$ Let $i+j=p.$ We need to
show that 
$$
H^1(M_C^{\otimes i}\otimes M_{rC}^{\otimes j}\otimes (p+k)C)=0 \ .
$$
We prove this by induction on $j.$ For $j=0,$ the vanishing is the
conclusion of Lemma 4.2. For $j>0,$ it is enough to prove that the following
multiplication map $\mu ,$ which sits in the following exact sequence

$$\displaylines{
H^0(M_C^{\otimes i}\otimes M_{rC}^{\otimes j-1}\otimes (p+k)C)\otimes
H^0(rC)@>{\mu }>> H^0(M_C^{\otimes i}\otimes
M_{rC}^{\otimes j-1}\otimes (p+r+k))C) \cr
\longrightarrow H^1(M_C^{\otimes i}\otimes M_{rC}^{\otimes j}\otimes
(p+k)C)\longrightarrow H^1(M_C^{\otimes i}\otimes M_{rC}^{\otimes
j-1}\otimes (p+k)C)\otimes H^0(rC),
}$$
surjects. The last term is zero by induction assumption, since $i+j-1=p-1.$
In view of Remark 2.1, it is enough to show that the following
multiplication map $\lambda $ sitting in the following exact sequence

$$\displaylines{
H^0(M_C^{\otimes i}\otimes M_{rC}^{\otimes j-1}\otimes (p+k)C)\otimes H^0(C)%
@>{\lambda }>>H^0(M_C^{\otimes i}\otimes
M_{rC}^{\otimes j-1}\otimes (p+k+1)C)
\cr
\longrightarrow H^1(M_C^{\otimes i+1}\otimes M_{rC}^{\otimes j-1}\otimes
(p+k)C),
}$$
surjects for all $k\geq 0$. This surjection follows from the vanishing of 
$$
H^1(M_C^{\otimes i+1}\otimes M_{rC}^{\otimes j-1}\otimes (p+k)C)
$$
which follows in turn by induction assumption on $j.$  $\square$ 

\proclaim{Lemma 4.4}
Let $X$ be a K3 surface and $L$ a line bundle on $X$. Assume that the
general member $C\in \mid L\mid $ is a smooth non-trigonal, non-plane
quintic curve of genus $g\geq 5.$ Then for any integer $p\geq 2,$

$$
H^1(M_C^{\otimes i}\otimes (i-1+k)C)=0\text{ for all }i=2,...,p\text{ and }%
k\geq 0.
$$
\endproclaim

{\bf {}\ Proof: }We will prove the lemma by induction on $i.$ If $i=2,$ we
want 
$$
H^1(M_C^{\otimes 2}\otimes (2-1+k)C)=0.
$$
This is true by Proposition 2.4. Assume the statement of the lemma for 
$p-1.$
So we have $H^1(M_C^{\otimes (p-1)}\otimes (p-2+k)C)=0$ for all $k\geq 0.$

We want to show that $H^1(M_C^{\otimes p}\otimes (p-1+k)C)=0.$ Tensoring
(2.1.1) with $M_C^{\otimes (p-1)}\otimes (p-1+k)C$ and taking global
sections yields the following sequence:
$$
\matrix
H^0(M_C^{\otimes (p-1)}\otimes (p-1+k)C)\otimes H^0(C) & \longrightarrow  & 
H^0(M_C^{\otimes (p-1)}\otimes (p+k)C) \\ 
\longrightarrow H^1(M_C^{\otimes p}\otimes (p-1+k)C) & \longrightarrow  & 
H^1(M_C^{\otimes (p-1)}\otimes (p-1+k)C)\otimes H^0(C)
\endmatrix
$$

The last term is zero by induction assumption. Since 
$$
H^1(M_C^{\otimes p-1}\otimes (p-2+k)C)=0
$$
by induction, we may use Lemma 2.1 to reduce the problem of seeing the
surjectivity of the map $\alpha $, which is a multiplication map of global
sections of vector bundles on a K3 surface, to the problem of checking the
surjectivity of a multiplication map of sections of bundles on a curve.
According to Lemma 2.1, it suffices to show that 
$$
H^0(M_C^{\otimes p-1}\otimes (p-1+k)C\otimes \Cal{O}_C)\otimes
H^0(C\otimes \Cal{O}_C)\longrightarrow H^0(M_C^{\otimes p-1}\otimes
(p+k)C\otimes \Cal{O}_C)
$$
surjects. This will follow from the vanishing 
$$
H^1(C,\text{ }M_C^{\otimes p-1}\otimes \Cal{O}_C\otimes M_{K_C}\otimes
(p-1+k)K_C)=0,
$$
where $C\otimes \Cal{O}_C=K_C$ is the canonical bundle on $C.$ To see
the vanishing, we use Lemma 3.1. The first condition required by Lemma 3.1
is the vanishing of $H^1($ $\Cal{O}_X).$ The second is the vanishing of $%
H^1((p-1+k)K_C),$ which occurs because $p\geq 3.$ Finally, to verify the
third condition one has to check the vanishing of 
$$
H^1(M_{K_C}^{\otimes p}\otimes (p-1+k)K_C)\ .
$$

Since the bundle $M_{C\otimes \Cal{O}_C}$ is stable (see [PR]) and the
tensor product of semistable bundles is semistable (See [Mi]), $M_{C\otimes 
\Cal{O}_C}^{\otimes p}\otimes (p-1)C\otimes \Cal{O}_C$ is semistable
so it is enough to check that $\mu (M_{C\otimes \Cal{O}_C}^{\otimes
(p)}\otimes (p-1+k)C\otimes \Cal{O}_C)>2g-2,$ where $g$ is the genus of $%
C$. We have\ 
$$
\mu (M_{C\otimes \Cal{O}_C}^{\otimes (p)}\otimes (p-1+k)C\otimes 
\Cal{O}_C)=p\mu (M_{K_C})+(p-1+k)\text{deg (}K_C)\ ,
$$
since $C\otimes \Cal{O}_C=K_C.$ Therefore we need $%
p(-2)+(p-2+k)(2g-2)>0, $ which is true for $p\geq 3$ and $g\geq 5.$\ 
$\square$

\medskip\ 

We are now ready to prove Theorem 4.1:

{\bf {}\ Proof }(of Theorem 4.1): We want to apply Lemma 3.3. For that we
have to check the two hypotheses in the statement of the lemma. The first
follows from Lemma 4.3 and the second follows from Lemma 4.4. \ $\square$
\medskip\ 

We now use the above results to prove more general vanishing theorems on K3
surfaces.\smallskip\ 

\proclaim{Lemma 4.5}
Let $X$ be a K3 surface, $L=\Cal{O}_X(C)$ a globally generated line
bundle such that the smooth general member $C\in \mid L\mid $ is a
non-hyperelliptic curve of genus $g>3.$ Then the cohomology groups 
$$\displaylines{
H^1(M_C^{\otimes j}\otimes M_{i_1C}^{\otimes j_1}\otimes M_{i_2C}^{\otimes
j_2}\otimes \cdots \otimes M_{i_rC}^{\otimes j_r}\otimes (p+k)C) = 0\text{
for all }i_s\geq 1\text{ } \cr
\text{and }j+j_1+j_2+\cdots +j_r = p\text{ , where }s=1,...,r.
}$$
\endproclaim

{\bf {}\ Proof: }We will prove the theorem only for $r=2$ (for the sake of
simplicity). The proof is very similar to the proof of Lemma 4.3. The proof
is by induction on $j+j_1+j_2$. If $j+j_1+j_2=1,$ then the theorem is true
by Proposition 2.2. Assume the theorem to be true for $p-1.$ Now let $%
j+j_1+j_2=p.$ So we want 
$$
H^1(M_C^{\otimes j}\otimes M_{i_1C}^{\otimes j_1}\otimes M_{i_2C}^{\otimes
j_2}\otimes (p+k))=0\ .
$$
We will prove the vanishing by induction on $j_1.$ If $j_1=0,$ the result
follows from Lemma 4.3. To prove for $j_1>0,$ it is enough to prove that $%
\mu ,$ which sits in the exact sequence 

$$\displaylines{
H^0(E\otimes (p+k)C)\otimes H^0(i_1C)@>{\mu }>>
H^0(E\otimes (p+k+i_1)C)
\cr
\longrightarrow H^1(F\otimes (p+k))\longrightarrow H^1(M_C^{\otimes
j}\otimes M_{i_1C}^{\otimes j_1-1}\otimes M_{i_2C}^{\otimes j_2}\otimes
(p+k)C)\otimes H^0(i_1C),
}$$
surjects, where $E=M_C^{\otimes j}\otimes M_{i_1C}^{\otimes j_1-1}\otimes
M_{i_2C}^{\otimes j_2}$ and $F=M_C^{\otimes j}\otimes M_{i_1C}^{\otimes
j_1}\otimes M_{i_2C}^{\otimes j_2}$. The last term is zero by induction
assumption, since $j+j_1-1+j_2=p-1.$

In view of Remark 2.1, it is enough to show the surjectivity of the
multiplication map $\lambda $ sitting in the following exact sequence: 
$$
H^0(E\otimes (p+k)C)\otimes H^0(C)@>{\lambda }>>
H^0(E\otimes (p+k+1)C)\longrightarrow H^1(G\otimes (p+k)C),
$$
where $G=$ $M_C^{\otimes j+1}\otimes M_{i_1C}^{\otimes j_1-1}\otimes
M_{i_2C}^{\otimes j_2}.$ The surjection follows from the vanishing of 
$$
H^1(M_C^{\otimes j+1}\otimes M_{i_1C}^{\otimes j_1-1}\otimes
M_{i_2C}^{\otimes j_2}\otimes (p+k)C),
$$
which in turn follows from induction assumption on $j_{1}$. \ $\square$

\medskip\ 

\proclaim{Theorem 4.6}
Let $X$ be a K3 surface and $L=\Cal{O}_X(C)$ be a line bundle such that
the general member $C\in \mid L\mid $ is a smooth non-trigonal, non-plane
quintic curve of genus $g\geq 5.$ Then 
$$\displaylines{
H^1(M_C^{\otimes j}\otimes M_{i_1C}^{\otimes j_1}\otimes M_{i_2C}^{\otimes
j_2}\otimes \cdots \otimes M_{i_rC}^{\otimes j_r}\otimes (p+k)C) = 0\text{
for all }j+j_1+\cdots +j_r=p+1 \cr
\text{and }i_t \geq 1\text{ where }t=1,...,r\text{ and }k\geq 0.
}$$
\endproclaim

{\bf {}Proof: }We will prove only the case $r=2,$ the general case is
exactly similar. So we have to show that 
$$
H^1(M_C^{\otimes j}\otimes M_{i_1C}^{\otimes j_1}\otimes M_{i_2C}^{\otimes
j_2}\otimes (p+k)C)=0\text{ for all }j+j_1+j_2=p+1,\text{ }i_s\geq 1\ .
$$
The proof follows from induction on $j_1+j_2.$ If $j_1+j_2=1,$ the result is
true by Theorem 4.1$.$ We assume the result for $j_1+j_2-1$ we want now to
prove it for $j_1+j_2.$ We need only to prove that the multiplication map in
the following long exact sequence surjects:
$$
\displaylines{
H^0(E\otimes (p+k)C)\otimes H^0(i_2C)  \longrightarrow   H^0(E\otimes
(p+k+i_1)C) \cr 
\longrightarrow H^1(F\otimes (p+k))  \longrightarrow   H^1(M_C^{\otimes
j}\otimes M_{i_2C}^{\otimes j_1}\otimes M_{i_2C}^{\otimes j_2-1}\otimes
(p+k)C)\otimes H^0(i_1C),
}
$$
where $E=M_C^{\otimes j}\otimes M_{i_1C}^{\otimes j_1}\otimes
M_{i_2C}^{\otimes j_2-1}$ and $F=M_C^{\otimes j}\otimes M_{i_1C}^{\otimes
j_1}\otimes M_{i_2C}^{\otimes j_2}.$ The last term in the above exact
sequence is zero by previous lemma, since $j+$ $j_1+j_2-1=p.$ In the light
of Remark 2.1 and arguments used repeatedly throughout this article it is
enough to prove that the multiplication map below is surjective: 
$$
H^0(E\otimes (p+k)C)\otimes H^0(C)\longrightarrow H^0(E\otimes
(p+k+1)C)\longrightarrow H^1(R\otimes (p+k)C)
$$
where $R$ $=$ $E\otimes M_C=M_C^{\otimes j+1}\otimes M_{i_1C}^{\otimes
j_1}\otimes M_{i_2C}^{\otimes j_2-1}.$ By induction assumption the result is
true for $j_1+j_2-1,$ hence we have the vanishing of 
$$
H^1(M_C^{\otimes j+1}\otimes M_{i_1C}^{\otimes j_1}\otimes M_{i_2C}^{\otimes
j_2-1}\otimes (p+k)C).\,\,\square 
$$
\medskip

\heading 5.Higher syzygies of K3 surfaces \endheading

In this section we give an application of the vanishing theorems proved in
Section 4. In particular, we show that the vanishing theorems imply a result
on higher syzygies for K3 surfaces.

Note that we are all along abusing the notation by writing $M_C$ instead of 
$%
M_L$, where $L=\Cal{O}_X(C).$ We now revert back to the notation $M_L$!
So for instance $M_{L^{\otimes r}}$ corresponds to the notation $M_{rC}$
which we have been using above.

\proclaim {Theorem 5.1}
Let $X$ be a K3 surface and $L=\Cal{O}_X(C)$ be a line bundle such that
the smooth general member $C\in \mid L\mid $ is a non-hyperelliptic curve of
genus $g>3.$ Then $L^{\otimes (p+1)}$ satisfies $N_p.$
\endproclaim

{\bf {}\ Proof: }Let $L^{\prime }=L^{\otimes (p+1)}.$ By Proposition 2.2 $%
L^{\prime }$ satisfies $N_0$. By [GL1], it is enough to prove 
$$
H^1({\bigwedge }^{i}M_{L^{\prime }}\otimes (L^{\prime })^{\otimes
s})=0\text{ for all }1\leq i\leq p+1\text{ and }s\geq 1.
$$
Since we are working over the complex numbers, it is enough to prove $%
H^1(M_{L^{\prime }}^{\otimes i}\otimes (L^{\prime })^{\otimes s})=0$ for all 
$1\leq i\leq p+1$ and $s\geq 1$. The vanishing follows from Lemma 4.3. \ $%
\square$ 

\proclaim{Theorem 5.2}
Let $X$ be a K3 surface and $L$ a line bundle on $X$. Assume further that
the general member $C\in \mid L\mid $ is a smooth non-trigonal, non-plane
quintic curve of genus $g\geq 5.$ Then $L^{\otimes p}$ satisfies $N_{p.}$
\endproclaim

{\bf {}\ Proof:} We need only to prove $H^1(M_{L^{\prime }}^{\otimes
i}\otimes (L^{^{\prime }})^{\otimes s})=0$, for all $1\leq i\leq p+1,$ where
now $L^{\prime }=L^{\otimes p}.$ By Lemma 4.3 and Theorem 4.1 we have the
vanishing of 
$$
H^1(M_{L^{\otimes r}}^{\otimes i}\otimes M_L^{\otimes j}\otimes L^{\otimes
(p+k)})
$$
for all $i+j\leq p+1$ and for all $k\geq 0.$ By letting $j=0$ and $r=p$, we
get the desired vanishing. \ $\square$ \medskip\ 

\heading  6. Syzygies of Fano varieties \endheading

In this section we prove a vanishing theorem on Fano threefolds of index
one, namely Theorem 6.6, and we generalize the result to Fano varieties of
dimension $n$ with index ($n-2)$. These results have as consequence results
on higher syzygies of Fano threefolds and, generally, of Fano varieties of
dimension $n$ with index $(n-2).$ The techniques and method of the proofs
are almost entirely analogous to the case of K3 surfaces, so most of the
time we will just sketch the proof and leave the details to the reader. %
\medskip\ 

Along the first part of this section (until Theorem 6.8), the variety $X$
will be a Fano threefold with very ample anticanonical bundle. We will write 
$L=-$ $K_X$ . Under these hypotheses a general member $S$ of $\mid L\mid $
is a smooth K3 surface. Note that if $X$ has index $1$, then $L$ is the
ample primitive line bundle generating Pic $(X)$.

For interesting examples of Fano 3-folds of index 1, we refer the reader to
[C].

\medskip

\proclaim{Lemma 6.1}
Let $X$ be a Fano threefold, $L=-$ $K_X$ very ample and $S$ a general member
of $\mid L\mid $. Then 
$$
H^1(M_{pS}\otimes rS)=0\text{ for all }r,\text{ }p\geq 1
$$
\endproclaim

{\bf {}\ Proof:} First we remark that, since $X$ is a Fano threefold, we
have $H^1(\Cal{O}_X(rS^{\prime }))=0$ for all $r\geq 1$, for any smooth
surface $S^{\prime }\subset X$. So tensoring the sequence 
$$
0\longrightarrow M_{pS}\longrightarrow H^0(pS)\otimes \Cal{O}%
_X\longrightarrow pS\longrightarrow 0
$$
by $rS$ yields: 
$$
H^0(pS)\otimes H^0(rS)@>{\varphi }>>%
H^0((p+r)S)\longrightarrow H^1(M_{pS}\otimes rS)\longrightarrow 0.
$$
Thus it is enough to prove that $\varphi $ is surjective. In view of Remark
2.1, it is enough to prove that the map 
$$
H^0(pS)\otimes H^0(S)\longrightarrow H^0((p+1)S)
$$
is surjective for all $p\geq 1$. Since $H^1(\Cal{O}_X)=0$ and $%
H^1(p^{\prime }S)=0$ for all $p^{\prime }\geq 1$, by Lemma 2.1 it is enough
to show that the map below is surjective: 
$$
H^0(pS\otimes \Cal{O}_S)\otimes H^0(S\otimes \Cal{O}%
_S)\longrightarrow H^0((p+1)S\otimes \Cal{O}_S)\ .
$$
This follows from Proposition 2.2 (Note that the general member $C\in \mid 
\Cal{O}_S(S)\mid $ is a non-hyperelliptic curve of genus $g\geq 3)$. \ $%
\square $ \medskip\ 

If we allow $p=1,$ we obtain Iskovskih's result (see [I]). \medskip

\proclaim{Remark 6.1}
Lemma 6.1 admits a slight variant in its statement: we could just assume
that $L$ is globally generated and that $C$ is non-hyperelliptic of genus $g$
$\geq 3.$ Under these assumptions 6.1 would imply in particular that $L$ is
very ample. All the results that follow can be reformulated as well in the
same fashion, thus indicating how our syzygy results on Fano varieties fit
in our meta-principle. 
\endproclaim
\medskip\
\proclaim{Lemma 6.2}
Let $X$ be a Fano threefold, $L=-K_X$ very ample and $S$ a general member of 
$\mid L\mid $ . Assume further that the smooth general member $C\in \mid 
\Cal{O}_S(S)\mid $ is a curve of genus $g>3.$ Then the cohomology group 
$$
H^1(X,\text{ }M_S^{\otimes p}\otimes (p+k)S)=0\text{ for all }k\geq 0,\text{ 
}p\geq 1.
$$
\endproclaim

{\bf {}\ Proof: }The proof as usual is by induction on $p.$ If $p=1$ the
vanishing holds by Lemma 6.1. Assume that the result is true for $p-1$; we
have 
$$
H^1(M_S^{\otimes p-1}\otimes (p-1+k)S)=0\text{ for all }k\geq 0.
$$
We want to show that $H^1(M_S^{\otimes p}\otimes (p+k)S)=0$. For that it is
enough to prove that the multiplication map $\alpha $ in the long exact
sequence

$$\displaylines{
H^0(M_S^{\otimes p-1}\otimes (p+k)S)\otimes H^0(S)@>{\alpha
}>> H^0(M_S^{\otimes p-1}\otimes (p+1+k)S) \cr
\longrightarrow
H^1(M_S^{\otimes p}\otimes (p+k)S)\longrightarrow H^1(M_S^{\otimes
p-1}\otimes (p+k)S)\otimes H^0(S),
}$$
surjects, since the last term of the sequence is zero by induction.

By Lemma 2.1 it is enough to show that $H^1(M_S^{\otimes p-1}\otimes
(p-1+k)S)=0$ in order to reduce the problem of checking the surjectivity of
the multiplication map $\alpha $ on the Fano threefold $X$ to the problem of
checking the surjectivity of a multiplication map on the K3 surface $S.$ The
required vanishing follows from induction hypothesis. So we need to prove
that the following multiplication map on $S$ surjects:

$$
H^0(M_S^{\otimes p-1}\otimes (p+k)S\otimes \Cal{O}_S)\otimes
H^0(S\otimes \Cal{O}_S)\longrightarrow H^0(M_S^{\otimes p-1}\otimes
(p+1+k)S\otimes \Cal{O}_S)
$$
the above map fits in the long exact sequence, whose next term is 
$$
H^1(M_S^{\otimes p-1}\otimes M_{S\otimes \Cal{O}_S}\otimes (p+k)S\otimes 
\Cal{O}_S).
$$
In view of Lemma 3.1 we need only to show that $H^1(S,$ $M_C^{\otimes
p}\otimes (p+k)C)$ $=0$. The vanishing follows from Lemma 4.2. $\square $%
\medskip\ 

\proclaim{Lemma 6.3}
Let $X$ be a Fano threefold, $L$ and $S$ as above with the additional
assumption that the smooth general member $C\in \mid \Cal{O}_S(S)\mid $
is a non-trigonal, non-plane quintic curve of genus $g\geq 5.$ Then

$$
H^1(X,M_S^{\otimes 2}\otimes kS)=0 \text{for all} k\geq 1.
$$
\endproclaim
{\bf {}\ Proof: }We will prove the Lemma for $k=1,$ the rest are similar
and easier. The cohomology group $H^1(M_S^{\otimes 2}\otimes S)$ is the
cokernel of the multiplication map

$$
H^0(M_S\otimes S)\otimes H^0(S)@>{\varphi }>>
H^0(M_S\otimes 2S),
$$
because $H^1(M_S\otimes S)=0$ by Lemma 6.1. Hence it is enough to prove that 
$\varphi $ surjects. Let $E=M_S\otimes S.$ Consider the commutative diagram:

$$
\matrix
0\longrightarrow & H^0(E)\otimes H^0(\Cal{O}_X) & \longrightarrow & 
H^0(E)\otimes H^0(S)\longrightarrow & H^0(E)\otimes H^0(S\otimes \Cal{O}%
_S)\longrightarrow 0 \\ 
& \downarrow &  & \downarrow & \downarrow \\ 
0\longrightarrow & H^0(E) & \longrightarrow & H^0(E\otimes S)\longrightarrow
& H^0(E\otimes S\otimes \Cal{O}_S)\longrightarrow 0
\endmatrix
$$

The left hand side of the above diagram surjects. In order to prove that the
multiplication map on the right hand side surjects, we use Lemma 2.1, since 
$%
H^1(\Cal{O}_X)=0$ and $H^1(M_S)=0$. Thus we reduce the problem of
checking the surjectivity of the multiplication map on the threefold to
checking on the K3 surface. By tensoring the sequence (see [GP2] ) 
$$
0\longrightarrow \Cal{O}_S\otimes H^0(\Cal{O}_X)\longrightarrow
M_S\otimes \Cal{O}_S\longrightarrow M_{(S\otimes \Cal{O}%
_S)}\longrightarrow 0
$$
by $S$ and taking global sections we have the following exact sequence:\ 
$$
0\longrightarrow H^0(S\otimes \Cal{O}_S)\longrightarrow H^0(M_S\otimes
S\otimes \Cal{O}_S)\longrightarrow H^0(M_{(S\otimes \Cal{O}%
_S)}\otimes S\otimes \Cal{O}_S)\longrightarrow 0.
$$
Tensoring the above sequence with $H^0(S\otimes \Cal{O}_S)$ and
considering the obvious multiplication maps yields the following commutative
diagram:

$$
\matrix
0\longrightarrow & H^0(C)\otimes H^0(C) & \longrightarrow & H^0(E\otimes 
\Cal{O}_S)\otimes H^0(C)\longrightarrow & H^0(F)\otimes
H^0(C)\longrightarrow 0 \\ 
& \downarrow &  & \downarrow & \downarrow \\ 
0\longrightarrow & H^0(2C) & \longrightarrow & H^0(E\otimes C\otimes 
\Cal{O}_S)\longrightarrow & H^0(F\otimes C)\longrightarrow 0
\endmatrix
$$
where $C\in \mid \Cal{O}_S(S)\mid $ and $F=M_{(S\otimes \Cal{O}%
_S)}\otimes S\otimes \Cal{O}_S.$ The maps on the left hand side and
right hand side are surjective by Propositions 2.2 and 2.4. $\square $ %
\medskip\ 

\proclaim{Lemma 6.4}
Let $X$ and $L=\Cal{O}_X(S)$ be as in Lemma 6.3. Then for any integer $%
p\geq 2,$ 
$$
H^1(M_S^{\otimes i}\otimes (i-1+k)S)=0\text{ for all }i=2,...,p\text{ and }%
k\geq 0.
$$
\endproclaim

{\bf {}\ Proof: }We will prove the lemma by induction on $i$ . If $i=2,$ we
want 
$$
H^1(M_S^{\otimes 2}\otimes (2-1+k)S)=0\text{ for all }k\geq 0.
$$
This is true by Lemma 6.3. Assume the statement of the lemma for $p-1.$ So
we have $H^1(M_S^{\otimes (p-1)}\otimes (p-2+k)S)=0$ .We want to show that 
$%
H^1(M_S^{\otimes p}\otimes (p-1+k)S)=0.$ The above group fits in the
following long exact sequence:

$$\displaylines{
H^0(M_S^{\otimes (p-1)}\otimes (p-1+k)S)\otimes H^0(S)@>{\beta
}>> H^0(M_S^{\otimes (p-1)}\otimes (p+k)S) \cr \longrightarrow
H^1(M_S^{\otimes p}\otimes (p-1+k)S)\longrightarrow H^1(M_S^{\otimes
(p-1)}\otimes (p-1+k)S)\otimes H^0(S).
}$$
Since the last term in the above sequence is zero by induction, it is enough
to show that $\beta $ is surjective. By Lemma 2.1 we can reduce the problem
of checking the surjectivity of $\beta $ resticting it to the K3 surface 
$S$%
, provided $H^1(M_S^{\otimes (p-1)}\otimes (p-2+k)S)=0,$ which is true again
by induction. So we want the following multiplication map on $S$ to surject:

$$
H^0(M_S^{\otimes (p-1)}\otimes (p-1+k)S\otimes \Cal{O}_S)\otimes
H^0(S\otimes \Cal{O}_S)\longrightarrow H^0(M_C^{\otimes (p-1)}\otimes
(p+k)C).
$$
The map fits in a long exact sequence whose next term is 
$$
H^1(M_S^{\otimes (p-1)}\otimes (p-1+k)S\otimes \Cal{O}_S\otimes
M_{S\otimes \Cal{O}_S})
$$
In order to show that the above cohomology group vanishes, it is enough, by
Lemma 3.1 to show that 
$$
H^1(S,\text{ }M_C^{\otimes p}\otimes (p-1+k)C)=0,.
$$
but this holds by Lemma 4.4. $\square $ \medskip\ 

\proclaim{Lemma 6.5}
Let $X$ be a Fano threefold, $L=-K_X$ very ample and $S$ a general member of 
$\mid L\mid $. Assume further that the smooth general member $C$ of $\mid 
\Cal{O}_S(S)\mid $ is a curve of genus $g>3.$ Then 
$$
H^1(M_S^{\otimes i}\otimes M_{rS}^{\otimes j}\otimes (p+k)S)=0\text{ for all 
}i+j=p,\text{ }k\geq 0\text{ and }r\geq 1.
$$
\endproclaim

{\bf {}\ Proof:} Mimic the proof of Lemma 4.3 word by word. $\square $%
\medskip\ 

\proclaim{Theorem 6.6}
Let $X$ be a Fano threefold, $L=-K_X$ very ample and $S$ a general member of 
$\mid L\mid $. Assume further that the smooth general member $C\in \mid 
\Cal{O}_S(S)\mid $ is a non-trigonal curve of genus $g\geq 5$, which is
not isomorphic to a smooth plane quintic. Then $H^1(M_S^{\otimes i}\otimes
M_{rS}^{\otimes j}\otimes (p+k)S)=0$ for all $i+j=p+1,$ $k\geq 0$ and $r\geq
1.$
\endproclaim

{\bf {}\ Proof: }We want to apply Lemma 3.3 for the Fano threefold $X$. The
hypothesis (3.3.1) and (3.3.2) needed in Lemma 3.3 follows from Lemmas 6.5
and 6.4 respectively.$\square $\medskip 

Now we can use the vanishing results obtained so far to prove results on the
syzygies of Fano threefolds as anounced:\ 

\proclaim{Theorem 6.7}
Let $X$ be a Fano threefold, $L=-K_X$ very ample and $S$ a general member of 
$\mid L\mid $. Assume further that the smooth general member $C\in \mid 
\Cal{O}_S(S)\mid $ is a curve of genus $g>3.$ Then $L^{\otimes (p+1)}$
satisfies $N_p.$
\endproclaim

{\bf {}\ Proof:} Denote $L^{\prime }=L^{\otimes p+1}.$ By Lemma 6.1 $%
L^{\prime }$ satisfies $N_0.$ Then by [GL1] it is enough to prove $H^1(%
{\bigwedge }^{i}M_{L^{^{\prime }}}\otimes (L^{^{\prime }})^{\otimes
s})=0$ for all $1\leq i\leq p+1$ and $s\geq 1.$ Since we are working over
the complex numbers, it is enough to prove $H^1(M_{L^{^{\prime }}}^{\otimes
i}\otimes (L^{^{\prime }})^{\otimes s})=0$ for all $1\leq i\leq p+1$. But
the needed vanishing follows from Lemma 6.5. $\square $ \medskip\ 

\proclaim{Theorem 6.8}
Let $X$ be a Fano threefold, $L=-K_X$ very ample and $S$ a general member of 
$\mid L\mid $. Assume further that the smooth general member $C\in \mid 
\Cal{O}_S(S)\mid $ is a non-trigonal curve of genus $g\geq 5$ which is
not isomorphic to a smooth plane quintic curve$.$ Then $L^{\otimes p}$
satisfies $N_p.$
\endproclaim

{\bf {}\ Proof: }Denote $L^{\prime }=L^{\otimes p}.$ By Lemma 6.1 $%
L^{\prime }$ satisfies $N_0.$ We have proved the vanishing 
$$
H^1(X,M_S^{\otimes i}\otimes M_{rS}^{\otimes j}\otimes (p+k)S)=0\text{ for
all }i+j\leq p+1,k\geq 0\text{ and }r\geq 1
$$
By letting $i=0,$ and $r=p,$ we have

$$
H^1(X,M_{pS}^{\otimes j}\otimes (p+k)S)=0\text{ for all }1\leq j\leq p+1.
$$
Then the theorem follows from the result in [GL1] mentioned in the previous
proof. $\square $\medskip\ 

As indicated in the introduction of this section the statements of the above
theorems and lemmas are in particular statements on the primitive ample
bundle of an index 1, Picard number 1 Fano threefold with very ample
anticanonical bundle. For that reason they hold more generally for the
primitive bundle $L=\Cal{O}_X(H)$ (and indeed, any ample bundle) of any
Fano variety of dimension $n$, index $n-2$ and Picard number $1$ for which $%
L $ is very ample. In such a situation $-K_X$ $=$ $(n-2)H.$ The basic
observation is this one:\medskip\ 

{\bf {}\ Observation :} Let $X$ be a Fano $n-$fold of index $(n-2)$ and $H$
a primitive, very ample line bundle on $X$. Then a smooth member in the
linear system of $\mid $ $H$ $\mid $ is a Fano $(n-1)$ fold of index 
$(n-3);$
let us call the smooth member also $H,$ then by adjunction

$K_H=(H-(n-2)H)\otimes \Cal{O}_H=-(n-3)H\otimes \Cal{O}_H.$ \medskip%
\ 

Therefore all the vanishing theorems and syzygy results which follow will be
obtained by induction on the dimension of the Fano variety. We will just
sketch the proofs and leave the details to the reader.\medskip\ 

\proclaim{Lemma 6.9}
Let $X$ be a Fano $n-$fold of index ($n-2)$ and Picard number $1$ having a
primitive, very ample line bundle $L=\Cal{O}_X(H)$. Assume that the
general $n-1$ section $H^{(n-1)}$ is a curve of genus $\geq 3.$ Then, 
$$
H^1(M_{pH}\otimes rH)=0\text{ for all }r,\text{ }p\geq 1
$$
\endproclaim

{\bf {}\ Proof:} Note first that on a Fano $n-$ fold, for any smooth
irreducible divisor $H^{\prime }$ $\subset X$, we have $H^1(\Cal{O}%
_X(H^{\prime }))=0.$ Then since $L$ is very ample, $H^{(n-1)}$ is a
non-hyperelliptic curve, for the general member of the linear system
corresponding to $\mid H^{(n-1)}\mid $ on the surface obtained by
intersecting $H$ $(n-2)$ times is a K3 surface. This allows us to use the
results proved for K3 surfaces as first step of the induction on the
dimension of $X$.

Tensoring the sequence 
$$
0\longrightarrow M_{pH}\longrightarrow H^0(pH)\otimes \Cal{O}%
_X\longrightarrow pH\longrightarrow 0
$$
by $rH$ yields:

$$
H^0(pH)\otimes H^0(rH)@>{\varphi }>>
H^0((p+r)H)\longrightarrow H^1(M_{pH}\otimes rH)\longrightarrow 0
$$
so it is enough to prove that $\varphi $ is surjective. In view of Remark
2.1, it is enough to prove that the map

$$
H^0(pH)\otimes H^0(H)^{\otimes r}@>{\alpha }>>
H^0((p+r)H)
$$
is surjective and by the same remark it is enough to show that 
$$
H^0(pH)\otimes H^0(H)@>{\beta }>> H^0((p+1)H)
$$
is surjective for any $p\geq 1$. By Lemma 2.1 it is enough to show that the
map 
$$
H^0(pH\otimes \Cal{O}_H)\otimes H^0(H\otimes \Cal{O}%
_H)\longrightarrow H^0((p+1)H\otimes \Cal{O}_H)
$$
surjects, since being $X$ a Fano $n-$ fold, $H^1(nH)=0$ for all $n\geq 0$ .
The above surjection follows from the induction assumption on the dimension
of the Fano variety. $\square $\medskip\ 

\proclaim{Lemma 6.10}
Let $X$ be a Fano $n-$fold of index ($n-2)$ and Picard number 1 having a
primitive and very ample bundle $L=\Cal{O}_X(H)$. Assume that the
general $n-1$ section $H^{(n-1)}$ is a curve of genus $g>3.$ Then the
cohomology group 
$$
H^1(X,\text{ }M_H^{\otimes p}\otimes (p+k)H)=0\text{ for all }k,p\geq 0.
$$
\endproclaim

{\bf {}\ Proof: }Mimic the proof of Lemma 6.2 word for word by replacing $S$
with $H.$ As in the last paragraph of Lemma 6.2 we need the vanishing of 
$$
H^1(S,M_S^{\otimes p-1}\otimes M_{S\otimes \Cal{O}_S}\otimes
(p+k)S\otimes \Cal{O}_S),
$$
here by Lemma 3.1 we need only to show that 
$$
H^1(H,\text{ }M_{H\otimes \Cal{O}_H}^{\otimes p-1}\otimes M_{H\otimes 
\Cal{O}_H}\otimes (p+k)H\otimes \Cal{O}_H)=0.
$$
The needed vanishing follows from the induction assumption on the 
dimension.$%
\square $\medskip\ 

\proclaim{Lemma 6.11}
Let $X$ and $L$ be as above with the additional assumption that $H\ ^{(n-1)}$
is a non-trigonal, non-plane quintic curve of genus $g\geq 5.$ Then 
$$
H^1(X,\text{ }M_H^{\otimes 2}\otimes H)=0.
$$
\endproclaim

{\bf {}\ Proof:} Mimic the proof of Lemma 6.3 and use the induction
assumption on the dimension. Note that if one member in the linear system of 
$\mid H\ ^{(n-1)}\mid $ is non-trigonal, the general member is also
non-trigonal.$\square $\medskip\ 

\proclaim{Lemma 6.12}
Let $X$ and $L$ be as in the above lemma. Then for any integer $p\geq 2,$ 
$$
H^1(M_H^{\otimes i}\otimes (i-1+k)H)=0\text{ for all }i=2,...,p\text{ and }%
k\geq 0.
$$
\endproclaim

{\bf {}\ Proof: }We may assume by induction the vanishing for Fano
varieties of dimension $(n-1)$ with index $(n-3).$ Mimic the proof of Lemma
6.4. All the needed vanishings for Fano varieties of one dimension less
follows from induction assumption.$\square $\medskip\ 

Note that by mimicking the arguments in Lemma 4.4 (as was done in Lemma
6.5) we obtain, 
$$
H^1(X,M_H^{\otimes i}\otimes M_{rH}^{\otimes j}\otimes (p+k)H)=0\text{ for
all }i+j=p,\text{ for all }k\geq 0,\text{ }r\geq 1\text{ }(*)
$$
Then by applying Lemma 3.3 we obtain, 
$$
H^1(X,M_H^{\otimes i}\otimes M_{rH}^{\otimes j}\otimes (p+k)H)=0\text{ for
all }i+j\leq p+1,\text{ for all }k\geq 0,\text{ }r\geq 1\text{ }(\bullet )
$$
\medskip\ 

\proclaim{Theorem 6.13}
Let $X$ be a Fano $n-$fold of index ($n-2)$ and $L=\Cal{O}_X(H)$ be the
primitive bundle. Assume that $L$ is very ample and $H^{(n-1)}$ is a smooth
non-hyperelliptic curve of genus $g>3.$ Then $L^{\otimes (p+1)}$ 
satisfies $%
N_p.$
\endproclaim

{\bf {}\ Proof: }By [GL1] it is enough to prove $H^1({\bigwedge }^{i}%
M_{L^{^{\prime }}}\otimes (L^{^{\prime }})^{\otimes s})=0$ for all $1\leq
i\leq p+1$ and $s\geq 1.$ Since we are working over the complex numbers, it
is enough to prove $H^1(M_{L^{^{\prime }}}^{\otimes i}\otimes (L^{^{\prime
}})^{\otimes s})=0$ for all $1\leq i\leq p+1$, where $L^{^{\prime
}}=L^{\otimes (p+1)}.$ The theorem follows from Lemma 6.10.$\square $
\medskip%
\ 

\proclaim{Theorem 6.14}
Let $X$ be a Fano $n-$fold of index ($n-2)$ and $L=\Cal{O}_X(H)$ be the
primitive bundle. Assume further that $H\ ^{(n-1)}$ is a non-trigonal,
non-plane quintic curve of genus $g\geq 5.$ Then $L^{\otimes p}$ satisfies $%
N_p.$
\endproclaim

{\bf {}\ Proof: }We have proved the vanishing $H^1(X,$ $M_H^{\otimes
i}\otimes M_{rH}^{\otimes j}\otimes (p+k)H)=0$ for all $i+j\leq p+1,$ for
all $k\geq 0$ and $r\geq 1.$ By letting $i=0,$ and $r=p,$ we have 
$$
H^1(X,M_{pH}^{\otimes j}\otimes (p+k)H)=0\text{ for all }1\leq j\leq p+1.
$$
So by [GL1] $L^{\otimes p}$satisfies $N_p.$ $\square $\bigskip\ 

\heading  References \endheading

\roster

\item"
[C]" S.D. Cutkosky,{\it \ On Fano 3-folds}, Manuscripta math. 64, 189-204
(1989)

\item"[EL]" L. Ein \& R. Lazarsfeld, {\it Koszul cohomology and Syzygies of
Projective
 varieties}, Inv Math {\bf {}\ 111} (1993), no1, 51-67.

\item"[GP1]" F. Gallego \& B.P. Purnaprajna, {\it Normal presentation on
elliptic ruled

surfaces,} (to appear in J of Algebra).

\item"[GP2]" \hbox{\leaders \hrule  \hskip .6 cm}\hskip .05 cm ,{\it  Higher
Syzygies of elliptic ruled surfaces}, (to appear in J of Algebra)

\item"[GP3]" \hbox{\leaders \hrule  \hskip .6 cm}\hskip .05 cm , {\it 
Syzygies of
Surfaces and Calabi-Yau 3-folds}, Preprint.

\item"[GL1]" M. Green \& R. Lazarsfeld, {\it Some results on the syzygies of
finite sets and
algebraic curves}, Compositio Math. {\bf {}\ 67 ( }1989) 301-314.

\item"[GL2]" \hbox{\leaders \hrule  \hskip .6 cm}\hskip .05 cm , {\it A 
simple
proof of Petri's theorem on canonical
curves,} Geometry Today (129-142), 1985, Birkhauser, Ed: Arbarello 
et al

\item"[I]" V.A. Iskovskih, {\it Fano 3-folds 1}, Izv.Akad.Nauk SSSR, Vol 11 
(1977)
no.3.

\item"[M]" A. Mayer, {\it Families of K3 surfaces, }Nagoya Mathematical 
Journal
(1972)

\item"[Mi]" Y. Miyaoka, {\it The chern class and Kodaira dimension of a 
minimal
variety,} Alg Geo-Sendai 1985, Adv. Studies in Pure Math., Vol 10,
449-476.

\item"[PR]" K. Paranjape \& S. Ramanan, {\it On the canonical ring of a 
curve},
Alg. Geo.
and Commutative Algebra in honor of Nagata, Vol 2, 503-516.

\item"[Pa-Pu]" G. Pareschi \& B. P.Purnaprajna, {\it\ Canonical ring of a 
curve
is Koszul:
A Simple proof,} (to appear in Illinois J. of Math.)

\item"[StD]" B. Saint-Donat, {\it Projective models of K3 surfaces}, 
Amer. J. of
Math. {\bf {}\ 96},
(1974) 602-639.
\endroster
\enddocument